\documentclass[12pt]{article}
\usepackage{amssymb,amsmath}
\begin{document}
\date{}
\title{Rigidity of noncompact complete Bach-flat manifolds}
\author{Seongtag Kim}
 \maketitle

\begin{center}
Department of Mathematics Education, Inha University,\\
Incheon,
402-751 Korea
\\ \textrm{e-mail:}  \texttt{stkim@inha.ac.kr}
\end{center}
\begin{abstract}

Let $(M,g)$ be a  noncompact complete Bach-flat manifold with
positive Yamabe constant.  We prove that $(M,g)$ is flat if  $(M,
g)$ has zero scalar curvature and sufficiently small  $L_{2}$
bound of curvature tensor. When $(M, g)$ has nonconstant scalar
curvature, we prove that $(M, g)$ is conformal to the flat space
if $(M, g)$ has sufficiently small $L_2$ bound of curvature tensor
and $L_{4/3}$ bound of scalar curvature.

\end{abstract}

\newtheorem{theorem}{Theorem}
\newtheorem{lemma}[theorem]{Lemma}
\newtheorem{cor}[theorem]{Corollary}
\def\bea{\begin{eqnarray}}
\def\eea{\end{eqnarray}}
\newcommand{\be}{\begin{equation}}
\newcommand{\ee}{\end{equation}}
\newcommand{\beq}{\begin{eqnarray}}
\newcommand{\eeq}{\end{eqnarray}}
\newcommand{\beqn}{\begin{eqnarray*}}
\newcommand{\eeqn}{\end{eqnarray*}}
\newcommand{\pa}{\partial}
\newcommand{\pxi}{ {\pa \over \pa x_i}}
\newcommand{\pxj}{ {\pa \over \pa x_j}}
\newcommand{\pxk}{ {\pa \over \pa x_k}}
\newcommand{\pxl}{ {\pa \over \pa x_l}}
\newcommand\nm{\nonumber}
\newcommand\ep{\epsilon}
\newcommand\ew{\epsilon'}
\newtheorem{rmk}{Remark}

\newcommand{\ai}[1]{\inf_{u\in C^{\infty}(M)}\frac{\int\limits_{#1}{|\nabla
u|^2+C_n S \, u^2 \; dV_g} }
 {(\int\limits_{#1}{u^{ 2n/( n-2) }\;  dV_g)^{(n-2)/n}  }  } }\par
\newcommand{\aI}[1]{\inf\limits_{u\in C_0^{\infty}(M-B_r)}  \frac {\int
\limits_{#1}{|\nabla u|^2+C_n S \, u^2 \; dV_g} }
 { (\int\limits_{#1}{u^{ 2n/(n-2)}\;  dV_g)^{(n-2)/n}  }    }}
\newcommand{\ab}[2]{  \frac{  \int\limits_{#1} |\nabla  {#2}|^2+C_n S
{#2}^2 dV_g} {     (\int\limits_{#1}{{#2}^{ 2n/(n-2)}\;
dV_g)^{(n-2)/n}} }}
%%%%%%%%%%%%%%%%%%%%%%%%%%%%%%%%%%%%%%%%%%%%%%%%%%%%%%%%%%%%%%%%%%%%%

\section{Introduction}

 Let  $(M,g)$  be a
 noncompact complete Riemannian $4$-manifold  with  scalar
 curvature $R$, Weyl curvature $W$,
 Ricci curvature $R_{ij}$ and curvature tensor $Riem$.
  A metric is Bach-flat
 if  it is a critical metric of the
functional \bea g \longrightarrow \int_M |W|^2 \,  dV_g.
\label{def1} \eea  Bach-flat condition is equivalent to the
vanishing of Bach tensor $B_{ij}$, which is defined by \bea
B_{ij}\equiv\nabla ^k \nabla ^l W_{kijl}+{1 \over 2}R^{k
l}W_{kijl} \label{e11} .\eea \noindent (see \cite{bes}). Important
examples of Bach manifolds are Einstein manifolds, self-dual
(anti-self-dual) manifolds, conformally flat manifolds and
K$\ddot{a}$ler surfaces with zero scalar curvature (see
\cite{der}).

Einstein metrics and Bach-flat metrics share many important
properties. When the curvature of a given Einstein metric $(M, g)$
is sufficiently close to that of the constant curvature space, in
$L_{n \over2} $ sense, it is known that $(M, g)$ is isometric to a
quotient of the constant curvature space \cite{an, sh, si,hv, it}.
In this paper, we study this rigidity phenomena for noncompact
complete Bach-flat manifolds. First, we study rigidity of metrics
with positive Yamabe constant,  zero scalar curvature and small
$L_2$ bound of curvature.  Next we look for rigidity of
nonconstant scalar curvature spaces with a conformal change of the
given metric.

 There are known rigidity of Bach-flat metrics. For a compact Bach-flat manifold $(M, g)$
with positive Yamabe constant, Chang, Ji and Yang \cite{ajp}
proved that there is only finite diffeomorphism class with an
$L_2$ bound  of Weyl tensor,  and $(M, g)$ is conformal to the
standard sphere if $L_2$ norm of Weyl tensor is small enough.
  For a noncompact complete Bach-flat manifold $(M, g)$
with positive Yamabe constant and zero scalar curvature, Tian and
Viaclovsky \cite{tv} proved that $(M, g)$ is almost locally
Euclidean of order $0$ with $L_2$ bounds of curvature, bounded
first Betti number and the uniform volume growth for any geodesic
ball.

For rigidity of Bach-flat manifolds, we use an elliptic estimation
on the Laplacian of curvature tensor. For this, we introduce  the
Yamabe constant on $(M, g)$.
 Let  $(M,g)$ be a Riemannian manifold of dimension $n
\ge 3$ with scalar curvature $R$. The Yamabe constant $Q(M, g)$ is
defined by
$$Q(M, g)\equiv \inf_{0\neq u\in C^{\infty}_0 (M)}  { { {    {{4(n-1)}\over {(n-2)}} \int_M |\nabla u|^2+R_g u^2 \,  dV_g}
                  \over {\left( \int_M |u|^{   {2n}/(n-2)} dV_g
                  \right)^{(n-2)/n}}} } . $$ \noindent
$Q(M, g)$ is conformally invariant and any locally conformally
flat manifold and manifolds with zero scalar curvature satisfy
$Q(M, g)>0$  (see \cite{del}).

%%%%%%%%%%%%%%%%%%%%%%%%%%%%%%%%%%%%%%%%%%%%%%%%%%%%%%%%%%%%%%%%%%%%%%%
\section{Bach-flat metric with constant scalar curvature}

In this section, we study noncompact complete Bach-flat manifolds
with nonnegative constant scalar curvature. First, we consider
Bach-flat manifolds whose $L_2$ curvature norm is small.  By an
elliptic estimation for the Laplacian of curvature tensor, we
have:
\begin{theorem}\label{th1} Let  $(M,g)$  be a
 noncompact complete Bach-flat Riemannian $4$-manifold with zero scalar curvature  and $Q(M, g)>0$.
 Then there exists a small number $c_0$  such that if $\int_M
   |Riem|^{2} \,  dV_g
 \le c_0$, then $(M,g)$  is flat,
 i.e., $Riem=0$, where $Riem$ is curvature tensor.

\end{theorem}
\noindent {\it Proof. \ } We need to prove that $|Riem|=0$. The
Laplacian of curvature tensor is \bea \Delta
R_{ijkl}&=&2(B_{ijkl}-B_{ijlk}-B_{iljk}+B_{ikjl})+\nabla_i\nabla_k
R_{jl}-\nabla_i\nabla_l
R_{jk}  \nonumber \\
{}&&-\nabla_j\nabla_k R_{il}+\nabla_j\nabla_l R_{ik}
+g^{pq}(R_{pjkl}R_{qi}+R_{ipkl}R_{qj}). \label{e30} \eea where
$B_{ijkl}=g^{pr}g^{qs}R_{piqj}R_{rksl}$  (see \cite{ham}).
  Multiplying $ R_{ijkl}$ on (\ref{e30}), \bea
&{}&R_{ijkl} \Delta R_{ijkl}\nm \\
&=&2(B_{ijkl}-B_{ijlk}-B_{iljk}+B_{ikjl})R_{ijkl}+(\nabla_i\nabla_k
R_{jl}-\nabla_i\nabla_l
R_{jk}) R_{ijkl}\nm \\
{}&&+(-\nabla_j\nabla_k R_{il}+\nabla_j\nabla_l R_{ik}) R_{ijkl}
+g^{pq}(R_{pjkl}R_{qi}+R_{ipkl}R_{qj}) R_{ijkl}.\label{q1} \eea To
simplify notations, we will work in an orthonormal frame. By the
Bianchi identity, \bea \nabla^i R_{ijkl}= \nabla_k R_{jl}-\nabla_l
R_{jk}. \label{bi2}\eea For a  smooth compact supported function
$\phi$ and small $\ep>0$, we integrate the second term in
(\ref{q1}) \bea &{}&\int_M \phi^2
(\nabla_i\nabla_k R_{jl}-\nabla_i\nabla_l R_{jk}) R_{ijkl} \,  dV_g \nm\\
&=& -\int_M \nabla_i \phi^2 (\nabla_k R_{jl}-\nabla_l R_{jk})
R_{ijkl}+\phi^2 (\nabla_k R_{jl}-\nabla_l R_{jk}) \nabla_i
R_{ijkl} \,  dV_g \nm\\ &=& -\int_M \nabla_i \phi^2 \nabla_t
R_{tjkl} R_{ijkl}+\phi^2
 |\nabla_i R_{ijkl}|^2 \,  dV_g \label{e4} \\
&\ge& -\int_M {1 \over \ep}|\nabla \phi|^2 | R_{ijkl}|^2+ (1+
\ep)\phi^2
 |\nabla_i R_{ijkl}|^2 \,  dV_g . \label{e5}
 \eea
 \noindent Using the same method,
\bea &{}& \int_M \phi^2
 (-\nabla_j\nabla_k R_{il}+\nabla_j\nabla_l R_{ik}) R_{ijkl} \,  dV_g\nm
 \\&\ge & -\int_M {1 \over \ep}|\nabla \phi|^2 | R_{ijkl}|^2+ (1+
\ep)\phi^2
 |\nabla_j R_{ijkl}|^2 \,  dV_g . \label{e6} \eea
 \noindent
  The first and fourth terms in (\ref{q1}) are contractions  of
 cubic terms of curvature tensor which can be bounded by
 $c|Riem|^3$ for a constant $c$. In this paper, we use $c$ and $c'$ to
 denote some positive constant, which can be varied.
  By the Kato inequality, $$ |\nabla Riem|^2\ge
|\nabla|Riem||^2$$ \noindent and  \bea &{}&-\int_M \phi^2
|Riem|\triangle |Riem| \,  dV_g
\\&=& -\int_M
 \phi^2\left( |\nabla Riem|^2-|\nabla|Riem||^2+ R_{ijkl}\triangle
R_{ijkl} \right) dV_g
\\&\le&
 \int_M 2 \left(
  { 1 \over \ep}|\nabla \phi|^2 | R_{ijkl}|^2+ (1+
\ep)\phi^2
 |\nabla_i R_{ijkl}|^2 \right) +c |Riem|^3 \phi^2 \,  dV_g. \label{q2}
 \eea
\noindent  For a general Riemannian $n$-manifold,
the following hold:\bea (\delta W)_{jkl} &=& \nabla^i W_{ijkl} \nm  \\
&=&{ {(n-3)} \over {(n-2)}} \left(\nabla_k R_{jl}-\nabla_l
R_{jk}-{1 \over 6} \nabla_k R g_{jl}+{1 \over 6} \nabla_l R
g_{jk}\right) \label{bh1} \eea and \bea |\nabla^i W_{ijkl}|^2=({
{n-3}\over {n-2}})^2\left(|\nabla^i R_{ijkl}|^2 -{1\over 6}|\nabla
R|^2 \right). \label{wh3}
 \eea   Let $E_{ij}$ be the traceless Ricci tensor, i.e
 $E_{ij}= R_{ij}-{1 \over 4} R g_{ij}$.
\noindent  Multiplying $\phi E_{ij}$ on Bach-flat equation
(\ref{e11}) (cf. \cite[(3. 24)]{cgy2}),
 \bea 0&=& \int_M \phi^2 E_{ij} \left(\nabla_k \nabla_\ell
W_{ik\ell j} - {1\over 2}
        W_{ikj\ell} E_{k\ell}\right) \,  dV_g \label{w15} \\
&=&  \int_M \phi^2\left(\big\vert \delta W\big\vert^2
         - {1\over 2} W_{ikj\ell} E_{k\ell} E_{ij} \right)-2 \phi  \nabla_k \phi
         \nabla_l  W_{iklj} E_{ij} \,  dV_g\nm \\
&=&  \int_M \phi^2\left(\big\vert \delta W\big\vert^2
         - {1\over 2} W_{ikj\ell} E_{k\ell} E_{ij} \right)-2 \phi  \nabla_k \phi
         \nabla_l  W_{iklj} E_{ij} \,  dV_g \label{e40} \\
&\ge&  \int_M (1-\epsilon_2) \phi^2 \big\vert \delta W\big\vert^2
         - {1\over 2}\phi^2 W_{ikj\ell} E_{k\ell} E_{ij} -{ 1\over \epsilon_2}  |\nabla
         \phi|^2
         | E_{ij} |^2 \,  dV_g \label{e37} ,
         \eea
\noindent where   $W_{ikj\ell} g_{kl}=0$ is used in (\ref{w15}).
Therefore, \bea
 \int_M  \phi^2 \big\vert \delta W\big\vert^2 \le
 (1-\epsilon_2)^{-1} \int_M
          {1\over 2}\phi^2 W_{ikj\ell} E_{k\ell} E_{ij} +{ 1\over \epsilon_2}  |\nabla
         \phi|^2
         | E_{ij} |^2 \,  dV_g \label{e137}
         \eea
\noindent and \bea &{}&-\int_M \phi^2 |Riem|\triangle |Riem| \,
dV_g \nm \\&\le&
 \int_M 2 \Big[
  { 1 \over \ep}|\nabla \phi|^2 | R_{ijkl}|^2+ {{4(1+
\ep)}\over{1-\ep_2}}\phi^2
 W_{ikj\ell} E_{k\ell} E_{ij} +{{8(1+
\ep)}\over{(1-\ep_2) \ep_2}}  |\nabla
         \phi|^2
         | E_{ij} |^2  \nm  \\&{}& \ + {1 \over 3}(1+
\ep)\phi^2 |\nabla R|^2 \Big ] +c |Riem|^3 \phi^2 \, dV_g
.\label{q5}
 \eea
Note that second term in (\ref{q5}) is also cubic term of
curvature tensor.  Now we can bound all cubic terms in the above
equation by $c |Riem|^3$,  and the first and third terms  by
$c|\nabla \phi|^2 |Riem|^2 $ for a suitable constant $c$. Next we
use the fact that scalar curvature is zero. For simplicity of
notations, we let $u=|Riem|$. Using the Yamabe constant $\Lambda_0
\equiv Q(M,g)$, \bea &{}&\Lambda_0 \left(\int_M (\phi u)^4\, dV_g
\right)^{1/2}\nm  \\ &\le& \int_M| u\nabla \phi + \phi \nabla u
|^2 \,  dV_g +{1\over 6} R u^2 \phi^2 \,  dV_g\label{e31}\\
&\le& \int_M u^2|\nabla \phi|^2+|\nabla u|^2 \phi^2 + 2 u \phi
\nabla \phi \cdot \nabla u +{ 1\over 6} R u^2 \phi^2 \ \,  dV_g  \nm  \\
&\le & \int_M
(c+1) |\nabla \phi|^2 u^2 +c u^3 \phi^2 \,  dV_g   \nm \\
  &\le & \int_M
 (c+1)|\nabla \phi|^2 u^2 \,  dV_g \nm   \\&{}& \  +c \left(\int_M (\phi u)^4
  \,  dV_g \right)^{1/2} \left(\int_M u^2 \,  dV_g\right)^{1/2}
  .\label{e34}\eea
  Since $\int_M |Riem|^2 \,  dV_g$ is sufficiently small, there exists a constant $c'$
  such that \bea
c' (\int_M (\phi u)^{4} \,  dV_g)^{1/2} \le  \int_M
  |\nabla \phi|^2 u^2 \,  dV_g . \label{c24}
  \eea
\noindent Now we choose
$\phi$ as \bea \phi=\begin {cases}1 & \text{on} \  B_t\\
0 & \text{on} \  M-B_{2t}
\\
|\nabla \phi|\le {2 \over t} & \text{on} \  B_{2t}-B_t \\
\end{cases} \eea with $0\le \phi \le 1$ and $B_t=\{ x\in M| d(x, x_0) \le t\}$ for some fixed $x_0 \in M$.
From (\ref{c24}) \bea c' \left(\int_M u^4 \phi^4 \,  dV_g
\right)^{1/2} &\le&
{4\over{t^2}}\int_{B(2t)-B(t)} u^2  \,  dV_g  \nm  \\
&\le& {4\over{t^2}} { \sqrt{3} \over 2}c_0.\eea By taking $t\to
\infty$, we have $u=0$. Therefore $(M,g)$ is flat.

\vskip 0.3 true cm

 Next we consider complete Bach-flat metric with positive constant scalar
 curvature. Using an elliptic estimation for traceless Ricci
 tensor, we prove

 \begin{theorem} \label{th4} Let  $(M,g)$  be a
 noncompact complete Riemannian $4$-manifold  with nonnegative constant  scalar
 curvature $R$, Weyl curvature $W$ and traceless Ricci curvature $E_{ij}$.
 Assume that $(M,g)$ is Bach-flat and $Q(M, g)>0$.
  Then there exists a small number $c_0$ such that if $\int_M |W|^2+
  |E_{ij}|^2\, dV_g
 \le c_0$, then $(M,g)$  is an Einstein
 manifold.

\end{theorem}

There is no noncompact complete Einstein manifold of positive
scalar curvature. By Theorem \ref{th4}, we have an obstruction for
the existence of a noncompact complete Bach-flat manifold.
\begin{theorem}\label{thob}
Let  $(M,g)$  be a noncompact complete Riemannian $4$-manifold
with positive constant  scalar
 curvature $R$, Weyl curvature $W$ and traceless Ricci curvature
 $E_{ij}$. Assume that $(M,g)$ is Bach-flat and $Q(M, g)>0$.
  Then, there exists a  positive number $c_1$ such that
  $\int_M |W|^2+ |E_{ij}|^2 \, dV_g\ge c_1$.
 \end{theorem}

\vskip 0.3 true cm \noindent {\it Proof} of Theorem \ref{th4}. \

Let $E_{ij}$ be the traceless Ricci tensor and $|E|=|E_{ij}|$.
Using Bianchi identity, Bach tensor can be expressed in the
following way  (cf. \cite[(1.18)]{cgy}), \bea B_{ij}&=&-{1 \over
2}\triangle E_{ij} +{ 1\over 6}
 \nabla_i \nabla_j R-{1 \over {24}}\triangle R
 g_{ij}-E^{kl}W_{ikjl}+E^k_i E_{jk} \nm  \\
 &{}&-{1 \over 4} |E|^2 g_{ij}+{1 \over 6}
 R E_{ij} \eea
  By the
Kato inequality, $ |\nabla E|^2\ge |\nabla|E||^2$ and $ {\rm tr}
E^3 \le {1 \over \sqrt{3}} |E|^3$, there exists a positive
constant $c$ satisfying the following equation for a Bach-flat
metric

\bea |E|\triangle |E|&=& |\nabla E|^2-|\nabla|E||^2-2
E^{kl}W_{ikjl}E^{ij}+ 2 {\rm tr} E^3 +{1 \over 3} R|E|^2 \\
&\ge&-2 E^{kl}W_{ikjl}E^{ij}+2 {\rm tr} E^3 +{1 \over 3} R|E|^2 \\
&\ge&-c |W||E|^2-{2 \over \sqrt{3}} |E|^3 +{1 \over 3} R|E|^2
 .\label{e1}\eea

\noindent Let $u=|E|$.  Multiplying a smooth compact supported
function $ \phi $ to (\ref{e1}) and  integrating on $M$,
$$\int_M \phi^2 |\nabla u|^2+ 2 \phi u \nabla \phi \cdot \nabla u \, dV_g
\le \int_M c |W| u^2 \phi^2 + {2 \over \sqrt{3}}u^3 \phi^2 - {1
\over 3} R u^2 \phi^2 \, dV_g. $$ \noindent Using the Yamabe
constant, \bea &{}& \Lambda_0 \left(\int_M (\phi u)^4 dv_g
\right)^{1/2} \nm  \\ &\le& \int_M| u \nabla  \phi + \phi \nabla \
u
|^2 dV_g +{1\over 6} R u^2 \phi^2 \nm   \\
&\le & \int_M u^2|\nabla \phi|^2+ c |W| u^2 \phi^2 + {2 \over
\sqrt{3}}u^3  \phi^2  -{1\over 6} R u^2 \phi^2dV_g. \label{e3}\eea

\noindent Note that the second term of (\ref{e3}) is bounded by
$$ c \left(\int_M|W|^2dV_g \right)^{1/2}   \left(\int_M u^4 \phi^4
\ dV_g \right)^{1/2} $$ \noindent and the third term is bounded by
$${2 \over \sqrt{3}} \left(\int_M u^4 \phi^4 dV_g
\right)^{1/2}\left(\int_M u^2 dV_g \right)^{1/2} . $$ \noindent

\noindent Assume that \bea c \left(\int_M |W|^2 dV_g \right)^{1/2}
+ {2 \over \sqrt{3}}\left( \int_M u^2 dV_g\right )^{1/2} \le
\Lambda_0 .\label{e444} \eea  Then, three terms in the right hand
side of (\ref{e3}) can be absorbed in the left hand side.
Therefore, there exists a constant $c'>0$ such that \bea c'
\left(\int_M u^4 \phi^4 dV_g \right)^{1/2} \le \int_M u^2 |\nabla
\phi|^2  dV_g \label{e445}\eea Now we choose
$\phi$ as \bea \phi=\begin {cases}1 & \text{on} \  B_t\\
0 & \text{on} \  M-B_{2t}\\
|\nabla \phi|\le {2 \over t} & \text{on} \  B_{2t}-B_t \\
\end{cases} \eea with $0\le \phi \le 1$.
From (\ref{e444}) and ( \ref{e445}) \bea c' \left(\int_M u^4
\phi^4 dV_g \right)^{1/2} &\le&
{4\over{t^2}}\int_{B(2t)-B(t)} u^2  dV_g \nm \\
&\le& {4\over{t^2}} { \sqrt{3} \over 2}\Lambda_0.\eea By taking
$t\to \infty$, we have $u=0$. Therefore $(M,g)$ is Einstein.

\section{Bach-flat metric with nonconstant scalar curvature }
In this section, we study noncompact complete Bach-flat metric
with nonconstant scalar curvature.
 We apply a result of the Yamabe problem
on noncompact manifold to study rigidity. For a given manifold
$(M, g)$, we find a conformal metric $\overline{g}=u^{4/(n-2)}g$
whose scalar curvature is zero. This is equivalent to find a
solution for the following partial differential equation \bea
-\Delta_g u+ {1\over 6} R_g u=0 .\label{yb} \eea

 The following existence of a
conformal metric with zero scalar curvature was proved by Kim
\cite{kim}.

\begin{theorem} \label{th2}
 Let $(M,g)$ be a noncompact complete Riemannian manifold  of dimension $n\ge 3$ with scalar curvature $R$.
 Assume that $Q(M,g)>0$ and $\int_M |R|^{2n/(n+2)}
+ |R|^{n/2} \, dV_g < \infty $. Then, there exists a conformal
metric $\overline{g}=u^{4/(n-2)}g$ whose scalar curvature is zero.
Moreover, $u$  satisfies the following:
 \bea \int_M | \nabla (u-1)|^2
 +|u-1|^{2n/(n-2)} \, dV_g < \infty \eea
and \bea &\int_M | \nabla (u-1)|^2 +|u-1|^{2n/(n-2)} \, dV_g
 \to 0  \ \text{as} \nm \\
 &\int_M |R|^{2n/(n+2)} + |R|^{n/2} \, dV_g \to 0. \eea
\end{theorem}

\vskip 0.5 true cm By Theorem \ref{th2} and an elliptic estimation
for solutions of (\ref{yb}), new metric $(M, \overline{g})$ in
 Theorem \ref{th2} is also complete (cf. \cite[ch. 8]{gt}).
By a standard elliptic estimation, $C^{2, \alpha}$ norm of $u-1$
is bounded by  $L_{n/2}$ and $L_{2n/(n+2)}$ norm of $R$.
Therefore, in dimension $4$, if $(M, g)$ has sufficiently small
$L_2$ bound of $|Riem|$ and  $L_{4/3}$ bound of $R$, there is a
conformal metric $\overline{g}$ with zero scalar curvature and
small $L_2$ norm of $|Riem_{\overline{g}}|$ with respect to
metric ${\overline{g}}$. Applying Theorem \ref {th1} to new metric
$(M, \overline{g})$, we have:

\begin{theorem}\label{th3} Let  $(M,g)$  be a
 noncompact complete Bach-flat Riemannian $4$-manifold with scalar curvature $R$ and $Q(M, g)>0$.
 Then there exists a small number $c_0$ such that if $\int_M
   |Riem|^{2}+|R|^{4/3} \,  dV_g
 \le c_0$, then $(M,g)$  is conformal to a flat space.

\end{theorem}

\begin{rmk} The constant $c_0$ in Theorem
\ref{th1}, \ref{th4},  \ref{th3} and $c_1$  in Theorem \ref{thob}
depend on $Q(M, g)$.
\end{rmk}

\section*{Acknowledgments} This work was supported by the Korea Research Foundation
Grant funded by the Korean Government (MOEHRD)
KRF-2005-202-C00036. This work was initiated when the author
visited Princeton University. The author thanks P. Yang for his
help and comments on this work and  Princeton University for
hospitality.

\end{document}